\documentclass[preprint,12pt]{elsarticle}

\usepackage{amssymb}
\usepackage{amsmath}

\usepackage{nicematrix}

\journal{Mathematics and Computers in Simulation}

\begin{document}

\begin{frontmatter}



\title{Meshing method to build a centrosymmetric matrix to solve partial differential equations on an irreducible domain including a planar symmetry} 


\author{T. Thuillier\corref{cor1}} 

\affiliation{organization={ Université Grenoble Alpes, CNRS, Grenoble INP, LPSC-IN2P3, 38000 Grenoble, France},
            addressline={}, 
            city={Grenoble},
            postcode={38026}, 
            state={},
            country={France}}
      
\cortext[cor1]{corresp. author: thomas.thuillier@lpsc.in2p3.fr}

\begin{abstract}
A general method to generate a centrosymmetric matrix associated with the solving of partial differential equation (PDE) on an irreducible domain by means of a linear equation system is proposed. The method applies to any PDE for which both the  domain to solve and the boundary condition (BC) type accept a planar symmetry, while no conditions are required on the BC values and the PDE right hand size function. It is applicable to finite element or finite difference method (FDM). It relies both on the specific construction of a mesh having a planar symmetry and a centrosymmetric numbering of the mesh nodes used to solve the PDE on the domain. The method is exemplified with a simple PDE using FDM. The method allows to reduce the numerical problem size to solve by a factor of two, decreasing as much the computing time and the need of computer memory.
\end{abstract}



\end{frontmatter}



\begin{figure}[t]

\centering
\includegraphics[width=1.\textwidth]{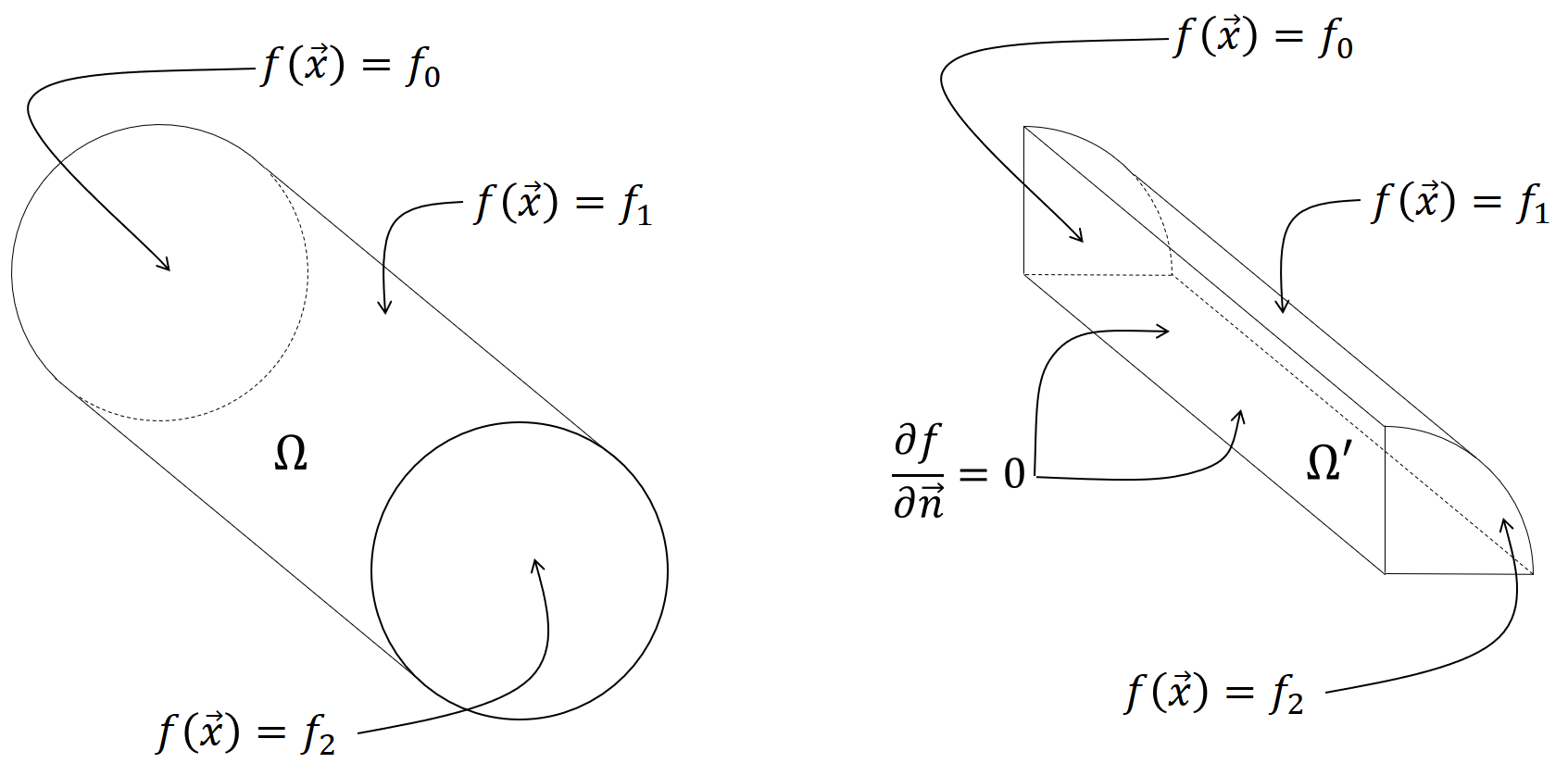}
\caption{Example of a domain reduction using known internal symmetries.}\label{fig:domain_1}

\end{figure}

\begin{figure}[t]

\centering
\includegraphics[width=0.5\textwidth]{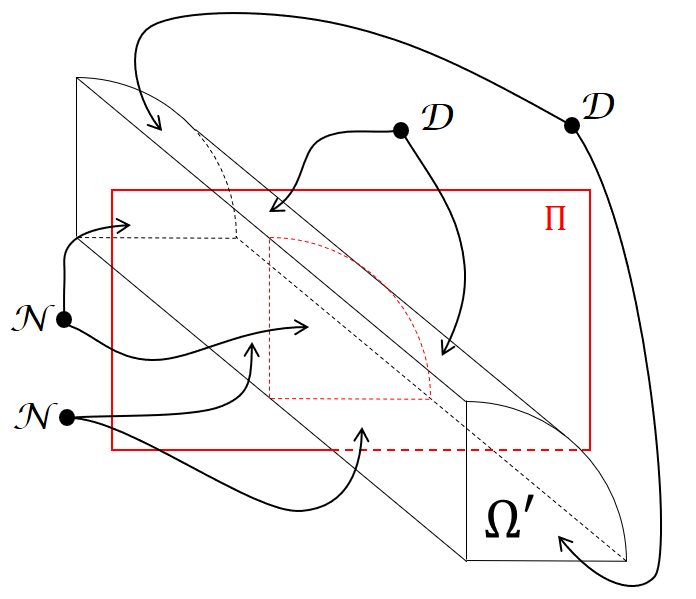}
\caption{Example of a reduced domain $\Omega'$ and boundary $\delta \Omega'$ with $\mathcal{N}$ and  $\mathcal{D}$ BC still presenting an unexploited planar geometrical symmetry (plane $\Pi$, red).}
\label{fig:domain_2}
\end{figure}
\section{Introduction}
\label{sec1}
 Physics and mathematics problems generally require the solving of partial differential equations (PDE)  on a domain $\Omega$ with specific boundary conditions on the border $\partial \Omega$~\citep{Evans1998}, having for instance the form
 \begin{equation}
     a_{ij}\frac{\partial^2 f}{\partial x_i \partial x_j}=g \, ,
     \label{eq:partial}
 \end{equation}
 where $a_{ij}$ and $g$ are known functions of $\vec{r}(x_1,x_2,x_3)$ and $f(\vec{r})$ is the function searched on $\Omega$. 
 Usual solving methods are the finite element (FEM)~\citep{Thomas2013} and the finite difference (FDM)~\citep{Ern2004}. They both rely on the discretization of space forming a mesh of $N$ points $M_n, n=1,..., N$  covering $\Omega$  and the linearization of the PDE to solve. The methods lead to the build up of a linear problem 
 \begin{equation}
 \label{eq:AXB}
     AX = b\, ,
 \end{equation}
 where $X$ is the unknown vector of component $f(M_n), n=1,..., N$, $b$ is a constant vector defined by specific values given on $\Omega$ and $\partial \Omega$ and $A$ is a square matrix of rank $N$. The  boundary conditions (BC) on $\partial \Omega$ are usually of Dirichlet type ($\mathcal{D}$), when $f$ is known, or Neumann ($\mathcal{N}$), when $\partial f / \partial \vec{n}$ is known ($\vec{n}$ normal to $\partial \Omega$ ), or a combination of the two.
When $\Omega$ is large or the need for accuracy is high, $N$  becomes rapidly very large ( $N \gg 10^6$) and the computer solving and the numerical storing of the inverse matrix $A^{-1}$, solution of Eq.~\ref{eq:AXB}, can become problematic. It is generally of interest to reduce the problem size by considering all the symmetries applying both on  $\Omega$, $\partial\Omega$ and $g$. An example of such a domain reduction from $\Omega$ down to $\Omega'$ is proposed in Fig.~\ref{fig:domain_1} for a specific cylindrical geometry, resulting in the onset of $\mathcal{N}$ BC. In many cases, the irreducible domain $\Omega'$ can still present a planar symmetry. Such a situation is exemplified in Fig.~\ref{fig:domain_2} where a plane of symmetry $\Pi$ is indicated in red. When in addition the BC types on $\Omega'$ admit the same symmetry, it is found that a centrosymmetric matrix $A$ can be created by carefully building the geometrical meshing and the mesh numbering. 

\section{Method to build a centrosymmetric matrix}
\label{sec2}

Let the total mesh number be $N'=2N$. Let $M_n$, $n=1,..., N$ be the ordered set of points of the mesh strictly located on one side of $\Pi$ ($M_n \cap \Pi =\emptyset$). Let $P_n, n=1,..., N$ be the image of the points $M_n$ by the plane $\Pi$. The complementary meshing located on the other side of $\Pi$ is then defined by the ordered set of points $M'_n$, with $n=N+1,..., 2N$, such that
\begin{equation}
\label{eq:centrosym}
M'_{N+1-n}=P_n , n=1,...,N\, ,
\end{equation}
leading to a centrosymmetric mesh numbering.
The linearization of Eq.~\ref{eq:partial} at any point $M_n$ creates the matrix row $A_n$ which geometrically relates $f_n=f(M_n)$ to its surrounding first neighbors $f_m=f(M_{m(n)})$. By construction, the matrix coefficients $a_{nm}$ involved in $A_n$ are only function of the mesh geometry (and this is true for both points in  $\Omega'$ and on $\delta \Omega'$). The equivalent relation obtained at the symmetrical point $P_n=M'_{N'+1-n}$ creates the matrix line $A_{N'+1-n}$. this line relates the function $f(P_n)$ to its values at the first neighbors $f(P_{m(n)})$. Because the  planar symmetry conserves the distances and the angles, the coefficients on the line $A_{N'+1-n}$ will have the same values as the one on the line $A_n$. And, by definition of Eq.~\ref{eq:centrosym}, the matrix element numbering involved on the row $A_{N'+1-n}$ will be centrosymmetric with respect to $A_n$: $a_{N'+1-n,N'+1-m}=a_{nm}$. 
It is worth noting that any geometry of mesh following the above conditions will result in the creation of a centrosymmetric matrix. In the case of FEM, the specific elements overlapping $\Pi$ must admit a planar symmetry with respect to $\Pi$ (like prisms with triangular of rectangular bases). \ref{ap:centrosym} proposes a few properties of centrosymmetric matrix and its inverse, showing that the computation time and the need of computer memory storage scales from $N'\,^2=4N^2$ down to $2N^2$, which is of high interest for the efficient solving of PDE with large $N'$.

\begin{figure}[t]

\centering
\includegraphics[width=0.7\textwidth]{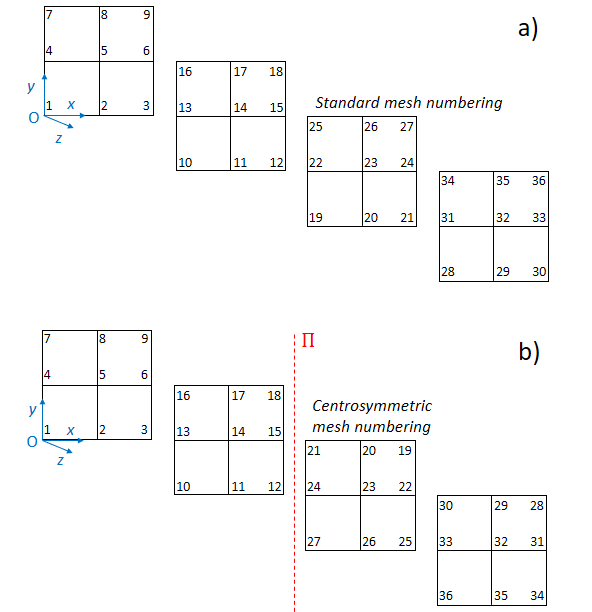}
\caption{a) Classical meshing numbering for the example treated. b) Centrosymmetric meshing numbering as proposed in this work. The plane of symmetry position (perpendicular to $z$) is indicated by a dashed red line.}\label{fig:meshing}
\end{figure}

\begin{figure}[t]
\centering
\includegraphics[width=0.7\textwidth]{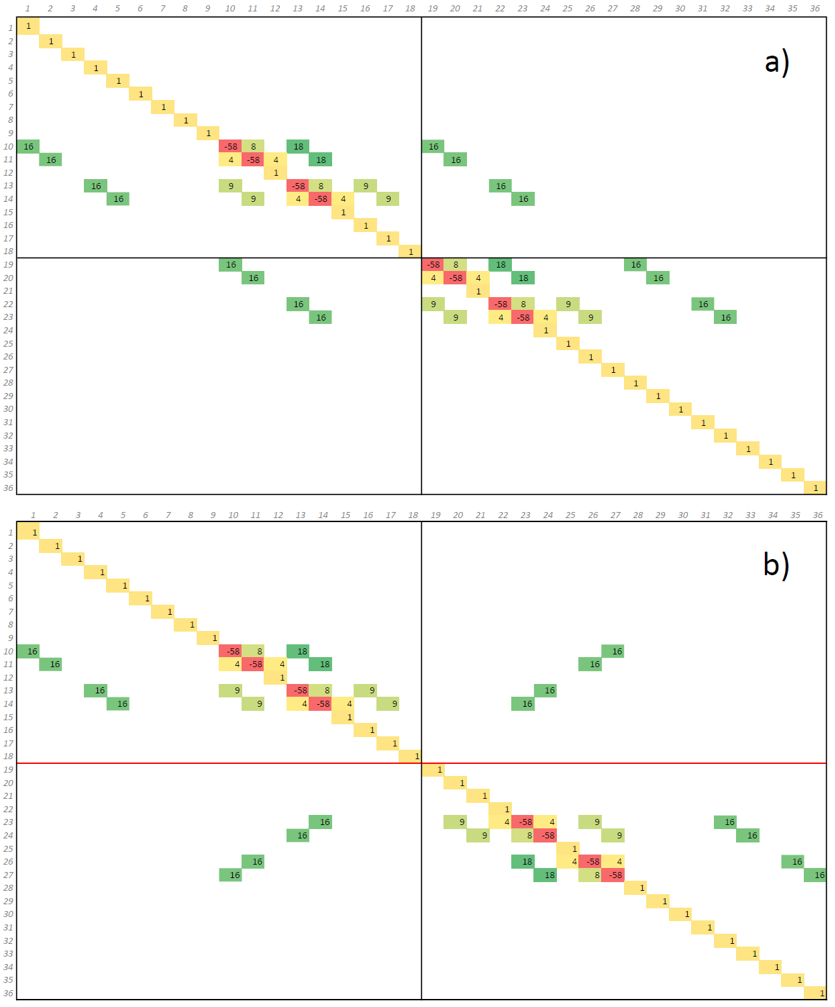}
\caption{a) Classical filling of matrix A. b) Centrosymmetric filling of matrix A. The matrix indexes are colored by value for visual convenience. The symmetry plane is indicated with a red line for b).}
\label{fig:matrix}
\end{figure}

\section{Example of matrix centro-symmetrisation}

This section proposes to calculate the matrix $A$ associated with the linearisation of equation 
\begin{equation}
    \label{eq:poisson}
\Delta f = \rho \, ,
\end{equation}
where $\rho$ is an arbitrary function. Both a classical and a centrosymmetric matrix $A$ are built for the sake of comprehension. $\Omega$ is chosen as a rectangular box fulfilling the planar symmetry condition. the Cartesian coordinate system $(O,x,y,z)$ is used. It is aligned with the three box orthogonal directions. $O$ is located on one edge of the box. The usual mesh numbering using a FDM is presented in Fig.~\ref{fig:meshing}(a) along with $O$ and the Cartesian vectors in blue. The mesh size along $x$, $y$ and $z$ are noted $h_x$, $h_y$ and $h_z$ respectively. A point $M_n$ of the mesh is identified by the three integers $(i, j, k)$ such that $\overrightarrow{OM}_n=i h_x\, \vec{x} + j h_y\, \vec{y} + k h_z\, \vec{z}$. The BC considered are of  $\mathcal{N}$ type for the left and bottom planes (respectively (xOy)  and (xOz) planes) and of $ \mathcal{D}$ type elsewhere. The Taylor expansion of Eq.\ref{eq:poisson} at $M_{n(i,j,k)}$ gives 
\begin{equation}
\label{eq:laplacian}
     \begin{split}
        h_{x}^{-2}\left(f_{i+1,j,k}+f_{i-1,j,k}\right)
        + h_{y}^{-2}\left(f_{i,j+1,k}+f_{i,j-1,k}\right)
        + h_{z}^{-2}\left(f_{i,j,k+1}+f_{i,j,k-1}\right)
        \\
        -2\left( h_{x}^{-2}+ h_{y}^{-2}+ h_{y}^{-2}\right) f_{i,j,k}
        = \rho_{i,j,k} \, ,
    \end{split}
\end{equation}
where $f_{i,j,k}=f(\overrightarrow{OM}_n)$ and $\rho_{i,j,k}=\rho(\overrightarrow{OM}_n)$. For $M_n$ on $\delta \Omega$, the $\mathcal{D}$ BC is expressed as 
\begin{equation}
    f_{ijk}=c_{ijk} \, ,
\end{equation}
where $c_{ijk}$ is an arbitrary local constant. The  $\mathcal{N}$ BC is expressed in $(xOZ)$ as 
\begin{equation}
\label{eq:neumann}
    f_{i,j+1,k}-f^*_{i,j-1,k}=2h_y\, q_{i,j,k} \, ,
\end{equation}
where $q_{i,j,k}$ is an arbitrary local flux,  $j=1$ and $f^*_{i,j-1,k}$ is a classical ghost function value out of $\Omega'$. Equation~\ref{eq:neumann} is used to eliminate the occurrence of $f^*_{i,j-1,k}$ in Eq.~\ref{eq:laplacian}. The $\mathcal{D}$ BC relation for the points located in the plane $(yOz)$ are accordingly expressed as 
\begin{equation}
    f_{i+1,j,k}-f^*_{i-1,j,k}=2h_x\, q'_{i,j,k} \, ,
\end{equation}
where $i=1$, $q'_{i,j,k}$ is an arbitrary local flux and $f^*_{i-1,j,k}$ a ghost function value treated as previously explained for the case $(xOz)$.
The equation system is built for the simple case when $i \in [0,2]$, $j \in [0,2]$ and $k \in [0,3]$, making $N'=36$ and $N=18$. $h_x$, $h_y$ and $h_z$ are chosen as $1/2$, $1/3$ and $1/4$ respectively to form integer matrix coefficients for convenience. The plane $\Pi$ is located at $z=\frac{3}{2}\,h_z$. With a classical mesh numbering (see Fig.~\ref{fig:meshing}.a), the matrix $A$ generated is plotted in Fig.~\ref{fig:matrix}.a. The typical band matrix structure of FDM is clearly visible and no symmetry is visible in the matrix. Now, when applying a centro-symmetrization of the meshing, as detailed in Fig.~\ref{fig:meshing}.b, the resulting interaction matrix $A$ is modified and the result is displayed in Fig.~\ref{fig:matrix}.b. One can see  that $A$ is now centrosymmetric. 

\section{Conclusions}
This work proposes a meshing method to create a centrosymmetric matrix associated with the linear solving of general PDE on an irreducible domain $\Omega$. the method is generally applicable provided that (i) $\Omega$ admits a planar symmetry with a plane $\Pi$, (ii) the BC types on $\delta \Omega$ are symmetric with $\Pi$ (while the local BC values can be different), (iii) the mesh network is symmetric with $\Pi$, (iv) The mesh numbering is centrosymmetric with respect to $\Pi$ and (v) no mesh points are located in $\Pi$. The method applies to FDM. The method applies to FEM provided that (vi) the elements overlapping $\Pi$ are symmetrical with $\Pi$ (e.g. prisms).  The matrix centro-symmetrization allows to reduce by a factor of two the computing time to inverse the matrix and the memory required to store the result in a computer memory.


\appendix
\section{Properties of centrosymmetric matrices}
\label{ap:centrosym}

A square centrosymmetric matrix with rank $2N$ can be expressed as a function of two block matrices $B$ and $C$ of rank $N$ ~\citep{GOOD1970}

\begin{equation}
A=\begin{pNiceArray}{c|c}
  B & CJ \\
  \hline
 JC & JBJ\\
\end{pNiceArray}\, ,
\end{equation}
where $J$ is the anti-diagonal identity matrix of rank $N$, with $J^2=I$ and $I$ the identity matrix of rank $N$. It is furthermore possible to demonstrate that $A^{-1}$ is centrosymmetric and expressed as

\begin{equation}
A^{-1}=\begin{pNiceArray}{c|c}
  D & EJ \\
  \hline
 JE &  JDJ\\
\end{pNiceArray}\, ,
\end{equation}
where $D$ and $E$ are matrices of rank $N$ with 
\begin{equation}
D=\frac{1}{2}\left((B+C)^{-1}+((B-C)^{-1}\right) \,
\end{equation}
and
\begin{equation}
E=\frac{1}{2}\left((B+C)^{-1}-((B-C)^{-1}\right).
\end{equation}

The matrix $A$ inversion of rank 2$N$ is thus obtained by inverting two matrices of rank $N$, leading to a computation cost scaling from $4N^2$ down to $2N^2$. Further information on centrosymmetric matrix inversion can be found  in~\citep{ElMikkawy2013,Khasanah2018}.

\bibliographystyle{elsarticle-num-names} 
\bibliography{bibliography.bib}

\begin{thebibliography}{6}
\expandafter\ifx\csname natexlab\endcsname\relax\def\natexlab#1{#1}\fi
\providecommand{\url}[1]{\texttt{#1}}
\providecommand{\href}[2]{#2}
\providecommand{\path}[1]{#1}
\providecommand{\DOIprefix}{doi:}
\providecommand{\ArXivprefix}{arXiv:}
\providecommand{\URLprefix}{URL: }
\providecommand{\Pubmedprefix}{pmid:}
\providecommand{\doi}[1]{\href{http://dx.doi.org/#1}{\path{#1}}}
\providecommand{\Pubmed}[1]{\href{pmid:#1}{\path{#1}}}
\providecommand{\bibinfo}[2]{#2}
\ifx\xfnm\relax \def\xfnm[#1]{\unskip,\space#1}\fi
\bibitem[{Evans(1998)}]{Evans1998}
\bibinfo{author}{L.~C. Evans},
\newblock \bibinfo{title}{Partial differential equations, second edition},
\newblock \bibinfo{journal}{Graduate Studies in Mathematics} \bibinfo{volume}{19} (\bibinfo{year}{1998}). \DOIprefix\doi{10.1090/gsm/019}.
\bibitem[{Ern and Guermond(2013)}]{Thomas2013}
\bibinfo{author}{A.~Ern}, \bibinfo{author}{J.~Guermond},
\newblock \bibinfo{title}{Theory and practice of finite elements.},
\newblock \bibinfo{journal}{Springer, Texts in Applied Mathematics} \bibinfo{volume}{22} (\bibinfo{year}{2013}). \DOIprefix\doi{/10.1007/978-1-4757-4355-5}.
\bibitem[{Evans(2004)}]{Ern2004}
\bibinfo{author}{L.~C. Evans},
\newblock \bibinfo{title}{Numerical partial differential equations: finite difference methods.},
\newblock \bibinfo{journal}{Springer, Applied Mathematical Sciences} \bibinfo{volume}{159} (\bibinfo{year}{2004}). \DOIprefix\doi{/10.1007/978-1-4899-7278-1}.
\bibitem[{Good(1970)}]{GOOD1970}
\bibinfo{author}{I.~J. Good},
\newblock \bibinfo{title}{The inverse of a centrosymmetric matrix},
\newblock \bibinfo{journal}{Technometrics} \bibinfo{volume}{12} (\bibinfo{year}{1970}) \bibinfo{pages}{925--928}. \DOIprefix\doi{10.1080/00401706.1970.10488743}.
\bibitem[{El-Mikkawy and Atlan(2013)}]{ElMikkawy2013}
\bibinfo{author}{M.~El-Mikkawy}, \bibinfo{author}{F.~Atlan},
\newblock \bibinfo{title}{On solving centrosymmetric linear systems.},
\newblock \bibinfo{journal}{Applied Mathematics} \bibinfo{volume}{4} (\bibinfo{year}{2013}) \bibinfo{pages}{21--32}. \DOIprefix\doi{10.4236/am.2013.412A003}.
\bibitem[{Khasanah et~al.(2018)Khasanah, Farikhin, and Surarso}]{Khasanah2018}
\bibinfo{author}{N.~Khasanah}, \bibinfo{author}{Farikhin}, \bibinfo{author}{B.~Surarso},
\newblock \bibinfo{title}{On computing: the inverse of centrosymmetric matrix},
\newblock \bibinfo{journal}{Journal of Physics: Conference Series} \bibinfo{volume}{983} (\bibinfo{year}{2018}) \bibinfo{pages}{012068}. \URLprefix \url{https://dx.doi.org/10.1088/1742-6596/983/1/012068}. \DOIprefix\doi{10.1088/1742-6596/983/1/012068}.

\end{thebibliography}
\end{document}